\input amstex 
\input xy
\xyoption{all}
\documentstyle{amsppt}

\document
\magnification=1200
\NoBlackBoxes
\nologo
\vsize18cm
\centerline{\bf Local action of the symmetric group}
\medskip 
\centerline{\bf and the twisted Yang-Baxter relation}
\bigskip
\centerline{\bf Alexander Odesskii}
\bigskip
\centerline{\bf Introduction}

\medskip

The Yang-Baxter relation plays a central role in the Quantum Field Theory. The usual $R$-matrix describes
the scattering of two "particles" (Fig.1). Here $u$ and $v$ are some characteristics of the
particles which are called the spectral parameters. $R$-matrix $(R_{ij}^{kl}(u,v))$ is a meromorphic function on the space of spectral parameters. Now suppose that the spectral parameters change after scattering (Fig.2). New spectral parameters $\varphi(u,v)$, $\psi(u,v)$ are some meromorphic functions in the variables $u$ and $v$. The factorization conditions impose a strong requirement on these functions (Fig.3). One can see from Fig.3 that

$$\varphi(u,\varphi(v,w))=\varphi(\varphi(u,v),\varphi(\psi(u,v),w))$$

$$\varphi(\psi(u,\varphi(v,w)),\psi(v,w))=\psi(\varphi(u,v),\varphi(\psi(u,v),w))$$

$$\psi(\psi(u,v),w)=\psi(\psi(u,\varphi(v,w)),\psi(v,w))$$

 If we have the solution of these equations, then we can define the analogous Yang-Baxter relation, $L$-operators, Zamolodchikov algebra and so on. In {\bf \S1} we develop the mathematical formalism for $\varphi$ and $\psi$ and give nontrivial examples. In {\bf \S2} we introduce the twisted Yang-Baxter relation and corresponding algebraic structures. The most interesting examples of $\varphi$ and $\psi$ appear in the theory of factorization of matrix polynomials ({\bf \S3}) and 
and matrix $\theta$-functions ({\bf \S4}). In {\bf \S5} we discus the generalized star-triangle relation and corresponding algebraic structures. In "physical" language it is the case when our particles are not elementary but are multiplets of $m$ elementary particles.

$$\xymatrix{v&u&\\u\ar[ur]&v\ar[ul]&} \xymatrix{&{\varphi(u,v)}&{\psi(u,v)}\\&u\ar[ur]&v\ar[ul]}$$
\phantom{aaaaaaaaaaaaaaaa} Figure 1 \phantom{aaaaaaaaaaaaaaaaaa} Figure 2

The examples of twisted $R$-matrices as well as the solutions of generalized star-triangle relation were found in [4] as intertwiners of cyclic representations and its tensor products of the algebra of monodromy matrices of the six-vertex model at root of unity [3]. These solutions are natural generalizations of the one from chiral Potts model [1,2,3]. One can obtain other solutions of the twisted Yang-Baxter and star-triangle relations by calculating the intertwiners of the representations of the algebras of monodromy matrices at root of unity for other trigonometric and elliptic $R$-matrices.

$$\xymatrix{&{\varphi(\varphi(u,v),\varphi(\psi(u,v),w))}&{\psi(\varphi(u,v),\varphi(\psi(u,v),w))}&{\psi(\psi(u,v),w)}\\ \\ \\ \\ \\ \\u\ar[uuuuuurrr]_(0.5){\psi(u,v)}&v\ar[uuuuuur]^(0.55){\varphi(u,v)}&&w\ar[uuuuuull]_(0.62){\varphi(\psi(u,v),w)}}$$
\medskip
\centerline{$\Vert$} 
\medskip
$$\xymatrix{{\varphi(u,\varphi(v,w))}&{\varphi(\psi(u,\varphi(v,w)),\psi(v,w))}&{\psi(\psi(u,\varphi(v,w)),\psi(v,w))}&\\ \\ \\ \\ \\ \\u\ar[uuuuuurr]^(0.62){\psi(u,\varphi(v,w))}&&v\ar[uuuuuul]_(0.55){\psi(v,w)}&w\ar[uuuuuulll]^(0.53){\varphi(v,w)}}$$

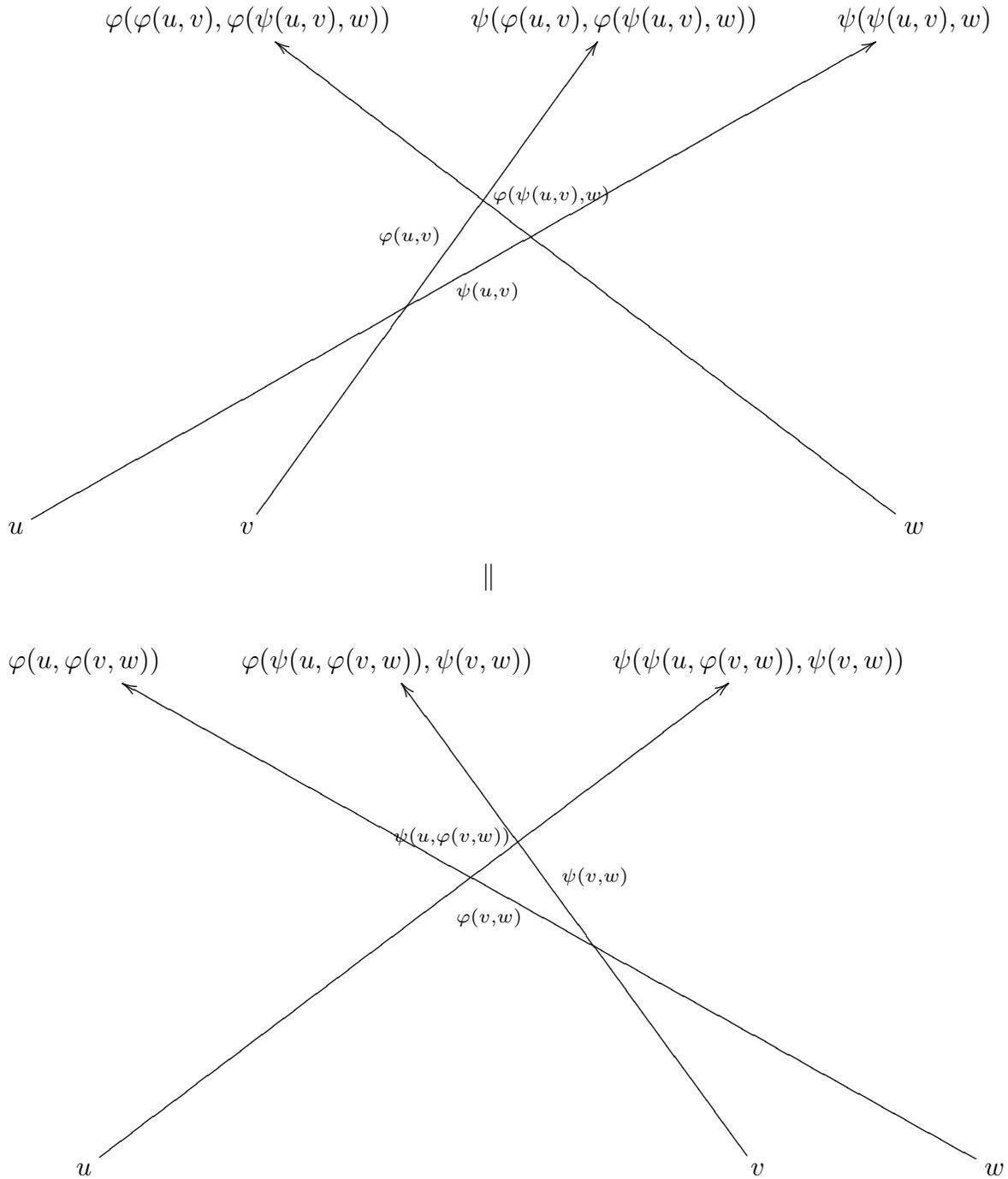
\centerline{Figure 3}

\newpage

\centerline{\bf \S1. Twisted transposition}
\medskip

Let $U$ be a complex manifold, $\mu: U\times U\to U\times U$ be a birational automorphism of $U\times U$. We will use the notations: $\mu(u,v)=(\varphi(u,v),\psi(u,v))$ where $u, v \in U$.
Here $\varphi$ and $\psi$ are meromorphic functions from $U\times U$ to $U$.

Let us introduce the following birational automorphisms of $U\times U\times U$: $\sigma_1=\mu\times \text{id}$ and $\sigma_2=\text{id}\times\mu$. We have: $\sigma_1(u,v,w)=(\varphi(u,v),\psi(u,v),w)$ and $\sigma_2(u,v,w)=(u,\varphi(v,w),\psi(v,w))$.

{\bf Definition} {\it We call a map $\mu$ a twisted transposition if the automorphisms $\sigma_1$ and $\sigma_2$ satisfy the following relations:}

$$\sigma_1^2=\sigma_2^2=\text{id},\text{     } \sigma_1\sigma_2\sigma_1=\sigma_2\sigma_1\sigma_2 \eqno(1)$$

If $\mu$ is a twisted transposition, then for each $N\in \Bbb N$ we have the birational action of the symmetric group $S_N$ on the manifold $U^N$ such that the transposition $(i,i+1)$ acts by automorphism $\sigma_i=\text{id}^{i-1}\times\mu\times \text{id}^{N-i-1}$. So we have $\sigma_i(u_1,\dots,u_N)=(u_1,\dots,\varphi(u_i,u_{i+1}),\psi(u_i,u_{i+1}),\dots,u_N)$. It is clear that the relations (1) are equivalent to the defining relations in the group $S_N$: $\sigma_i^2=e, \sigma_i\sigma_{i+1}\sigma_i=\sigma_{i+1}\sigma_i\sigma_{i+1}, \sigma_i\sigma_j=\sigma_j\sigma_i$ for $|i-j|>1$.

It is easy to check that the relations (1) are equal to the following functional equations for $\varphi$ and $\psi$:

$$\varphi(\varphi(u,v),\psi(u,v))=u,\text{   } \psi(\varphi(u,v),\psi(u,v))=v$$

$$\varphi(u,\varphi(v,w))=\varphi(\varphi(u,v),\varphi(\psi(u,v),w))$$

$$\varphi(\psi(u,\varphi(v,w)),\psi(v,w))=\psi(\varphi(u,v),\varphi(\psi(u,v),w)) \eqno(2)$$

$$\psi(\psi(u,v),w)=\psi(\psi(u,\varphi(v,w)),\psi(v,w))$$

{\bf Remarks 1.} From (2) it follows that for each $N$ the functions $\varphi(u_1,\varphi(u_2,\dots,\varphi(u_N,w)\dots)$ and $\psi(\dots(\psi(w,u_1),u_2)\dots,u_N)$ are invariant with respect to the action of the group $S_N$ on the variables $u_1,\dots,u_N$.

{\bf 2.} Informally one can consider $\mu$ as an infinite dimensional $R$-matrix in the space of functions. Namely, if we consider the space of meromorphic functions $\{f: U\times U\to \Bbb C\}$ as an "extended tensor square" of the space of meromorphic functions $\{f: U\to\Bbb C\}$, then the linear operator $R_{\mu}: f\to f\circ\mu$ (that is $R_{\mu}f(u,v)=f(\mu(u,v))$) satisfies the usual Yang-Baxter relation.

{\bf Examples 1.} Let $q, q^{-1}: U\to U$ be birational automorphisms such that $qq^{-1}=q^{-1}q=\text{id}$. Then 
$$\mu(u,v)=(q(v),q^{-1}(u))$$
is a twisted transposition.

{\bf 2.} Let $U=\Bbb C$, then the following formula gives a twisted transposition:

$$\mu(u,v)=(1-u+uv, \frac{uv}{1-u+uv})$$

{\bf 3.} Let $U$ be a finite dimensional associative algebra with unity $1\in U$ and such that generic element $u\in U$ is invertible, for example $U=\text{Mat}_m$. Then the following formula gives a twisted transposition:

$$\mu(u,v)=(1-u+uv, (1-u+uv)^{-1}uv)$$

For further examples see {\bf \S3,4}.

\newpage

\centerline{\bf \S2. Twisted Yang-Baxter relation and}
\centerline{\bf corresponding algebraic structures}
\medskip

Let $V$ be a $n$-dimensional vector space. For each $u\in U$ we denote by $V(u)$ a vector space canonically
isomorphic to $V$. Let $R$ be a meromorphic function from $U\times U$ to $End(V\otimes V)$. We will consider $R(u,v)$ as a linear operator 
$$R(u,v): V(u)\otimes V(v)\to V(\varphi(u,v))\otimes V(\psi(u,v))$$

{\bf Definition} {\it We call $R$ a twisted $R$-matrix (with respect to the twisted transposition $\mu$) if it satisfies the following properties:

1. The composition

$$V(u)\otimes V(v)\to V(\varphi(u,v))\otimes V(\psi(u,v))\to V(u)\otimes V(v)$$
is equal to the identity, that is $R(\varphi(u,v),\psi(u,v))R(u,v)=1$.

2. The following diagram is commutative:
$$\xymatrix{{V(\varphi(u,v))\otimes V(\psi(u,v))\otimes V(w)}\ar[r]^(.42){1\otimes R(\psi(u,v),w)}&{V(\varphi(u,v))\otimes V(\varphi(\psi(u,v)),w))\otimes V(\psi(\psi(u,v),w))}\ar[d]^{R(\varphi(u,v),\varphi(\psi(u,v),w))\otimes 1}\\{V(u)\otimes V(v)\otimes V(w)}\ar[u]_{R(u,v)\otimes 1}\ar[d]^{1\otimes R(v,w)}&{\widetilde V}\\{V(u)\otimes V(\varphi(v,w))\otimes V(\psi(v,w))}\ar[r]^(.42){R(u,\varphi(v,w))\otimes 1}&{V(\varphi(u,\varphi(v,w)))\otimes V(\psi(u,\varphi(v,w)))\otimes V(\psi(v,w))}\ar[u]_{1\otimes R(\psi(u,\varphi(v,w)),\psi(v,w))}}$$

Here $\widetilde V=V(\varphi(u,\varphi(v,w)))\otimes V(\psi(\varphi(u,v),\varphi(\psi(u,v),w)))\otimes V(\psi(\psi(u,v),w))$.

In other words, 

$$R^{12}(\varphi(u,v),\varphi(\psi(u,v),w)))R^{23}(\psi(u,v),w)R^{12}(u,v)=$$
$$R^{23}(\psi(u,\varphi(v,w)),\psi(v,w))R^{12}(u,\varphi(v,w))R^{23}(v,w) \eqno(3)$$

Here $R^{12}=R\otimes 1$ and $R^{23}=1\otimes R$ are linear operators in $V\otimes V\otimes V$.

We call (3) twisted Yang-Baxter relation.}

Let $\{x_i, i=1\dots,n\}$ be a basis of the linear space $V$, $\{x_i(u)\}$ be the corresponding basis of the linear space $V(u)$. It is clear that the following two linear operators are twisted $R$-matrices for each $\mu$:
$$x_i(u)\otimes x_j(v)\to x_i(\varphi(u,v))\otimes x_j(\psi(u,v))$$ 
$$x_i(u)\otimes x_j(v)\to x_j(\varphi(u,v))\otimes x_i(\psi(u,v))$$ 

Let $R^{\alpha\beta}_{ij}(u,v)$ be a matrix element of twisted $R$-matrix $R$, which means

$$R(u,v): x_i(u)\otimes x_j(v)\to R^{\alpha\beta}_{ij}(u,v)x_{\beta}(\varphi(u,v))\otimes x_{\alpha}(\psi(u,v))$$
We sum over repeating indices.

{\bf Definition} {\it We call the twisted Zamolodchikov algebra associated with $R$ an associative algebra $Z_R$ with generators $\{x_i(u); i=1\dots,n,u\in U\}$ and defining relations}

 $$x_i(u)x_j(v)=R^{\alpha\beta}_{ij}(u,v)x_{\beta}(\varphi(u,v))x_{\alpha}(\psi(u,v))$$

It is clear that the twisted Zamolodchikov algebra satisfies the usual property: for generic $u_1,\dots,u_N\in U$ the elements $\{x_{i_1}(u_1)\dots x_{i_N}(u_N); 1\leqslant i_1,\dots,i_N\leqslant n\}$ are linear independent.

Let $W$ be a vector space, $L$ be a meromorphic function $L: U\to End(V\otimes W)$. We will consider $L(u)$ as a linear operator $L(u): V(u)\otimes W\to W\otimes V(u)$. 

{\bf Definition} {\it We call $L$ an $L$-operator for twisted $R$-matrix $R$ if the following diagram is commutative:
$$\xymatrix{{V(u)\otimes W \otimes V(v)}\ar[rr]^{L(u)\otimes 1}&&{W\otimes V(u)\otimes V(v)}\ar[d]^{1\otimes R(u,v)}\\{V(u)\otimes V(v)\otimes W}\ar[u]_{1\otimes L(v)}\ar[d]^{R(u,v)\otimes 1}&&{W\otimes V(\varphi(u,v))\otimes V(\psi(u,v))}\\{V(\varphi(u,v))\otimes V(\psi(u,v))\otimes W}\ar[rr]^{1\otimes L(\psi(u,v))}&&{V(\varphi(u,v))\otimes W \otimes V(\psi(u,v))}\ar[u]_{L(\varphi(u,v))\otimes 1}}$$ 

Which is

$$R(u,v)L^1(u)L^2(v)=L^1(\varphi(u,v))L^2(\psi(u,v))R(u,v) \eqno(4)$$
where $L^1=L\otimes 1$ and $L^2=1\otimes L$.}

Let $\{L^j_i(u); i,j=1,\dots,n\}$ be the elements in $End(W)$ such that $L(u): x_i(u)\otimes W\to (L^{\alpha}_i(u)W)\otimes x_{\alpha}(u)$. We can consider $L(u)$ as a $n\times n$ matrix with matrix elements $(L^j_i(u))$. 

{\bf Definition} {\it The algebra of monodromy matrices $M_R$ for the twisted $R$-matrix $R$ is an associative algebra with generators  $\{L^j_i(u); i,j=1,\dots,n, u\in U\}$ and defining relations (4).}

It is clear that representations of the algebra $M_R$ in a vector space $W$ coincide with $L$-operators $L(u)\in End(V\otimes W)$.

One can introduce a coproduct $\Delta$ in the algebra $M_R$ by the usual way:

$$\Delta L^j_i(u)=L_i^{\alpha}(u)\otimes L^j_{\alpha}(u)$$
In the language of $L$-operators, if we have two $L$-operators $L_1(u): V(u)\otimes W_1\to W_1\otimes V(u)$ and $L_2(u): V(u)\otimes W_2\to W_2\otimes V(u)$, then we can construct an $L$-operator in the space $W_1\otimes W_2$ as the composition:

$$V(u)\otimes W_1\otimes W_2\to W_1\otimes V(u)\otimes W_2\to W_1\otimes W_2\otimes V(u)$$

{\bf Definition} {\it $Q$-operator is the element of the algebra $M_R$ given by the formula:}

$$Q(u)=\text{tr}L(u)=L^i_i(u)$$

From relations (4) it follows that $Q$-operators satisfy the following relations:

$$Q(u)Q(v)=Q(\varphi(u,v))Q(\psi(u,v))$$

{\bf Remark} One can also define the dynamical twisted Yang-Baxter relation.

\newpage
\centerline{\bf \S3. Local action of the symmetric group}
\centerline{\bf and factorization of matrix polynomials}
\medskip

For the general theory of matrix polynomials and factorization see [7]. For our purposes we state the results, which may be well known to the experts.

We denote by $S(a)$ the set of eigenvalues of a matrix $a\in Mat_m$. More generally, we denote by $S(a_1,\dots,a_d)$, $a_1,\dots,a_d\in Mat_m$, the set of roots of the polynomial $f(t)=\det(t^d-a_1t^{d-1}+\dots+(-1)^da_d)$. We will consider polynomials with generic coefficients only, so $\#S(a_1,\dots,a_d)=md$.

{\bf Proposition 1.} {\it Let 
$$t^d-a_1t^{d-1}+\dots+(-1)^da_d=(t-b_1)\dots(t-b_d) \eqno(5)$$
for generic matrices $a_1,\dots,a_d\in Mat_m$, then $S(b_i)\cap S(b_j)=\emptyset$ for $i\ne j$ and $S(b_1)\cup\dots\cup S(b_d)=S(a_1,\dots,a_d)$. For each decomposition $S(a_1,\dots,a_d)=A_1\cup\dots\cup A_d$, such that $\#A_i=m$, $A_i\cap A_j=\emptyset$ $(i\ne j)$ there exists a unique the factorization (5) with $S(b_i)=A_i$.}

{\bf Proof} The first statement follows from the equation $\det(t^d-a_1t^{d-1}+\dots+(-1)^da_d)=\det(t-b_1)\dots \det(t-b_d)$.

On the other hand, if we know eigenvalues of $b_1,\dots,b_d$ then we can calculate eigenvectors of them. For $\lambda\in S(b_d)$ the corresponding eigenvector is a vector $v_{\lambda}$, such that $({\lambda}^d-a_1{\lambda}^{d-1}+\dots+(-1)^da_d)v_{\lambda}=0$. If we know all eigenvectors of $b_d$, then we can calculate all eigenvectors of $b_{d-1}$ similarly and so on. This implies the uniqueness. By our construction of $b_1,\dots,b_d$ the determinants of the matrix polynomials in the right hand side and the left hand side of (5) have the same set of roots. Moreover, for each root $\lambda$ the operators represented by these matrix polynomials have the same kernel if we set $t=\lambda$. It implies that these polynomials are equal.

{\bf Proposition 2.} {\it Let $a_1, a_2\in Mat_m$ be generic matrices. Then there exists a unique pair of matrices $b_1, b_2\in Mat_m$ such that $(t-a_1)(t-a_2)=(t-b_1)(t-b_2)$ and $S(b_1)=S(a_2), S(b_2)=S(a_1)$. We have $b_1=a_1+{\Lambda}^{-1}, b_2=a_2-{\Lambda}^{-1}$ where $a_2\Lambda-\Lambda a_1=1$.}

{\bf Proof} If $S(b_2)\cap S(a_2)\ne\emptyset$, then $\det(a_2-b_2)=0$, because $a_2$ and $b_2$ have a common eigenvector. Otherwise, we can put $\Lambda=(a_2-b_2)^{-1}$.

Let $U=Mat_m$. From the propositions 1 and 2 it follows that the formula $\mu(a_1,a_2)=(b_1,b_2)$ gives a twisted transposition, where $b_1+b_2=a_1+a_2$, $b_1b_2=a_1a_2, S(b_1)=S(a_2), S(b_2)=S(a_1)$. We have $\mu(a_1,a_2)=(a_1+{\Lambda}^{-1},a_2-{\Lambda}^{-1})$, where $\Lambda$ is the solution of the linear matrix equation  $a_2\Lambda-\Lambda a_1=1$.

Let $\overline U=\overline{Mat_m}$ be the set of $m\times m$ matrices with different eigenvalues and fixed order of eigenvalues. The proposition 1 gives an action of the symmetric group $S_{mN}$ on the space $\overline U^N$ by birational automorphisms. By definition, for $\sigma\in S_{mN}$, $b_1,\dots,b_N\in \overline U$ we have $\sigma(b_1,\dots,b_N)=(b_1^\prime,\dots,b_N^\prime)$ where $(t-b_1)\dots(t-b_N)=(t-b_1^\prime)\dots(t-b_N^\prime)$ and $\overline S(b_i^\prime)=\sigma\overline S(b_i)$,  $\overline S$ stands for the ordered set of eigenvalues.

This action is local in the following sense. The transposition $(i,i+1)$ for $\alpha m<i<(\alpha+1)m$ acts only inside the $\alpha+1$-th factor of $\bar U^N$ and the transposition $(\alpha m,\alpha m+1)$ acts only inside the product of the $\alpha$-th and the $\alpha+1$-th factors. We have also the twisted transposition $\mu: \overline U\times \overline U\to \overline U\times \overline U$ in this case which is the action of the element $(1,m+1)(2,m+2)\dots(m-1,2m-1)\in S_{2m}$.

{\bf Remark} Let $U$ be the set of matrix polynomials of the form $at+b$, where $a,b\in Mat_m, a=(a_{ij}), b=(b_{ij})$ and $a_{ij}=0$ for $i<j$, $b_{ij}=0$ for $i>j$. It is possible to define the twisted transposition $\mu$ such that for $\mu(f(t),g(t))=(f_1(t),g_1(t))$ we have $f(t)g(t)=f_1(t)g_1(t)$, $\det f(t)$ and $\det g_1(t)$ have the same set of roots and the first coefficients of $f(t)$ and $g_1(t)$ have the same diagonal elements. In [4] we found the solutions of the corresponding twisted Yang-Baxter relation (for $m=2$), which is generalization of the $R$-matrix from chiral Potts model.

\newpage
\centerline{\bf \S4. Local action of the symmetric group }
\centerline{\bf and factorization on matrix $\theta$-functions}
\medskip

Let $\Gamma\subset\Bbb C$ be the lattice generated by 1 and $\tau$ where $\text{Im}\tau>0$. We have $\Gamma=\{\alpha+\beta\tau; \alpha,\beta\in\Bbb Z\}$. Let $\varepsilon\in\Bbb C$ be primitive root of unity of degree $m$. Let $\gamma_1,\gamma_2\in Mat_m$ be $m\times m$ matrices such that $\gamma_1^m=\gamma_2^m=1, \gamma_2\gamma_1=\varepsilon\gamma_1\gamma_2$. We have $\gamma_1v_{\alpha}=\varepsilon^{\alpha}v_{\alpha}, \gamma_2v_{\alpha}=v_{\alpha+1}$ in some basis $\{v_{\alpha}; \alpha\in\Bbb Z/m\Bbb Z\}$ of $\Bbb C^m$. Let us assume that $\{v_1,\dots,v_m\}$ is the standard basis of $\Bbb C^m$. 

We denote by $M\Theta_{n,m,c}(\Gamma)$ for $n,m\in\Bbb N, c\in\Bbb C$ the space of everywhere holomorphic functions $f: \Bbb C\to Mat_m$, which satisfy the following equations:

$$f(z+\frac{1}{m})=\gamma_1^{-1}f(z)\gamma_1$$
$$f(z+\frac{1}{m}\tau)=e^{-2\pi i(mnz-c)}\gamma_2^{-1}f(z)\gamma_2 \eqno(6)$$

{\bf Proposition 3.} {\it $\dim M\Theta_{n,m,c}(\Gamma)=m^2n$ and for each element $f\in M\Theta_{n,m,c}(\Gamma)$ the equation $\det f(z)=0$ has exactly $mn$ zeros modulo $\frac{1}{m}\Gamma$. The sum of these zeros is equal to $mc+\frac{mn}{2}$ modulo $\Gamma$.}

{\bf Proof} For $m=1$ we have the usual $\theta$-functions $\Theta_{n,c}(\Gamma)=M\Theta_{n,1,c}(\Gamma)$ and all these statements are well known in this case ([8]). One have a basis $\{\theta_{\alpha}(z); \alpha\in\Bbb Z/n\Bbb Z\}$ in the space  $\Theta_{n,c}(\Gamma)$ such that $\theta_{\alpha}(z+\frac{1}{n})=e^{2\pi i\frac{\alpha}{n}}\theta_{\alpha}(z), \theta_{\alpha}(z+\frac{1}{n}\tau)=e^{-2\pi i(z-\frac{n-1}{2n}\tau-\frac{1}{n}c)}\theta_{\alpha+1}(z)$ ([8]). From (6) it follows that $f(z+1)=f(z)$ and $f(z+\tau)=e^{-2\pi i(m^2nz-c_1)}f(z)$ for some $c_1\in\Bbb C$. So the matrix elements of $f(z)$ are $\theta$-functions from the space $\Theta_{m^2n,c_1}(\Gamma)$. We have decomposition $f(z)=\sum_{\alpha}\varphi_{\alpha}\theta_{\alpha}(z)$, where $\varphi_{\alpha}\in Mat_m$ are constant matrices, $\{\theta_{\alpha}\}$ is a basis in the space $\Theta_{m^2n,c_1}(\Gamma)$. Substituting this decomposition in (6) one can calculate the dimension of the space $M\Theta_{n,m,c}(\Gamma)$. We have also $\det f(z+\frac{1}{m})=\det f(z)$ and $\det f(z+\frac{1}{m}\tau)=e^{-2\pi i(m^2nz-mc)}\det f(z)$. From this follows the statement about zeros of the equation $\det f(z)=0$.

{\bf Proposition 4.} {\it For generic complex numbers $\lambda_1,\dots,\lambda_{mn}$ such that $\lambda_1+\dots+\lambda_{mn}\equiv mc+\frac{mn}{2}$ $\text{mod}\frac{1}{m}\Gamma$ and nonzero vectors $v_1,\dots,v_{mn}\in\Bbb C^m$ there exist a unique up to proportionality element $f(z)\in M\Theta_{n,m,c}(\Gamma)$ such that $\det f(\lambda_{\alpha})=0, f(\lambda_{\alpha})v_{\alpha}=0$ for $1\leqslant\alpha\leqslant mn$.}

{\bf Proof} Considering the decomposition $f(z)=\sum_{\alpha}\varphi_{\alpha}\theta_{\alpha}(z)$ one has the system of linear equations $\{\sum_{\alpha}\theta_{\alpha}(\lambda_{\beta})\varphi_{\alpha}v_{\beta}=0; \beta=1,\dots,mn\}$ for matrix elements of $\{\varphi_{\alpha}\}$. One can see that this system defines $\{\varphi_{\alpha}\}$ uniquely up to proportionality for generic $\lambda_1,\dots,\lambda_{mn}, v_1,\dots,v_{mn}$.

We denote by $S(f)$ the set of zeros of the equation $\det f(z)=0$ modulo $\frac{1}{m}\Gamma$.

{\bf Proposition 5.} {\it Let $f(z)\in M\Theta_{n,m,c}(\Gamma)$ be generic element and we have factorization $f(z)=f_1(z)\dots f_n(z)$ where $f_{\alpha}(z)\in M\Theta_{1,m,c_{\alpha}}(\Gamma)$, $c_1+\dots+c_n=c$. Then $S(f_{\alpha})\cap S(f_{\beta})=\emptyset$ for $\alpha\ne\beta$ and $S(f)=S(f_1)\cup\dots\cup S(f_n)$. For each decomposition $S(f)=A_1\cup\dots\cup A_n$ such that $A_{\alpha}\cap A_{\beta}=\emptyset$ for $\alpha\ne\beta$ and $\#A_{\alpha}=m$ there exists a unique factorization $f(z)=f_1(z)\dots f_n(z)$ up to proportionality of $f_{\alpha}$ such that $S(f_{\alpha})=A_{\alpha}$ for $\alpha=1,\dots,n$.}

{\bf Proof} is similar to the proof of proposition 1, we just change polynomials by $\theta$-functions.

Let $U_c$ be the projectivisation of the linear space  $M\Theta_{1,m,c}(\Gamma)$ and $U=\bigcup_{c\in\Bbb C}U_c$. We have the following twisted transposition $\mu: U\times U\to U\times U$. By definition $\mu(f,g)=(f_1,g_1)$, where $f(z)g(z)=f_1(z)g_1(z)$ and $S(f_1)=S(g), S(g_1)=S(f)$. 

Let $\overline U$ be the set of elements $f$ from $U$ with fixed order on $S(f)$. For $f\in\overline U$ let $\overline S(f)$ be the set $S(f)$ with corresponding order. We have a local action of the symmetric group $S_{mN}$ on the space $\overline U^N$. By definition, for $\sigma\in S_{mN}$ we have $\sigma(f_1,\dots,f_N)=(f_1^{\sigma},\dots,f_N^{\sigma})$ where $f_1(z)\dots f_N(z)=f_1^{\sigma}(z)\dots f_N^{\sigma}(z)$ and $\overline S(f_{\alpha}^{\sigma})=\sigma\overline S(f_{\alpha})$ for $1\leqslant\alpha\leqslant N$.

{\bf Remark} It is possible to construct twisted $R$-matrices for this twisted transposition $\mu$ as intertwiners of tensor products of cyclic representations of the algebra of monodromy matrices for elliptic Belavin $R$-matrix [5] at the point of finite order (see also [6]). It will be the subject of another paper. 

\newpage
\centerline{\bf \S5. Generalization of the star-triangle relation}
\medskip

Let $U$ be a complex manifold with an action of the symmetric group $S_m$ by birational automorphisms. Let $g_i$ be a birational automorphism of $U$ corresponding to the transposition $(i,i+1)$, where $1\leqslant i<m$. We assume also that the group $S_{2m}$ acts on the manifold $U\times U$ by birational automorphisms, such that the transposition $(i,i+1)$ for $i<m$ acts by $g_i\times \text{id}$ and the transposition $(i,i+1)$ for $i>m$ acts by $\text{id}\times g_{i-m}$. Let $f$ be the birational automorphism of $U\times U$ corresponding to the transposition $(m,m+1)$. We will use the notations $f(u,v)=(\lambda(u,v),\mu(u,v))$. We have: $g_i^2=f^2=\text{id}$, $g_ig_{i+1}g_i=g_{i+1}g_ig_{i+1}$, $fg_m^lf=g_m^lfg_m^l$, $fg_1^rf=g_1^rfg_1^r$, where $g_m^l=g_m\times \text{id}$ and $g_1^r=\text{id}\times g_1$.

It is clear that we have the local action of the symmetric group $S_{mN}$ on the manifold $U^N$ for each $N$.

Let $V$ be a $n$-dimensional vector space. For each $u\in U$ we denote by $V(u)$ a vector space canonically isomorphic to $V$. Let $G_i$ $(1\leqslant i<m)$ be a meromorphic function from $U$ to $End(V)$ and $F$ be a meromorphic function from $U\times U$ to $End(V\otimes V)$. We will consider $G_i(u)$ and $F(u,v)$ as linear operators 

$$G_i(u): V(u)\to V(g_i(u))$$
$$F(u,v): V(u)\otimes V(v)\to V(\lambda(u,v))\otimes V(\mu(u,v))$$

{\bf Definition} {\it We call $G_i, F$ a twisted $GF$-matrices if the following properties hold:

1. The compositions
$$V(u)\to V(g_i(u))\to V(u)$$
$$V(u)\otimes V(v)\to V(\lambda(u,v))\otimes V(\mu(u,v))\to V(u)\otimes V(v)$$
are equal to the identity.

2. The following diagrams are commutative:}
$$\xymatrix{{V(g_i(u))}\ar[rr]^{G_{i+1}(g_i(u))}&&{V(g_{i+1}g_i(u))}\ar[d]^{G_i(g_{i+1}g_i(u))}\\{V(u)}\ar[u]_{G_i(u)}\ar[d]^{G_{i+1}(u)}&&{V(g_ig_{i+1}g_i(u))}\\{V(g_{i+1}(u))}\ar[rr]^{G_i(g_{i+1}(u))}&&{V(g_ig_{i+1}(u))}\ar[u]_{G_{i+1}(g_ig_{i+1}(u))}}$$

$$\xymatrix{{V(\lambda(u,v))\otimes V(\mu(u,v))}\ar[r]^{G_m(\lambda(u,v))\otimes 1}&{V(g_m(\lambda(u,v)))\otimes V(\mu(u,v))}\ar[d]^{F(g_m(\lambda(u,v))),\mu(u,v))}\\{V(u)\otimes V(v)}\ar[u]_{F(u,v)}\ar[d]^{G_m(u)\otimes 1}&{V(\lambda(g_m(\lambda(u,v)),\mu(u,v)))\otimes V(\mu(g_m(\lambda(u,v)),\mu(u,v)))}\\{V(g_m(u))\otimes V(v)}\ar[r]^(.4){F(g_m(u),v)}&{V(\lambda(g_m(u),v))\otimes V(\mu(g_m(u),v))}\ar[u]_{G_m(\lambda(g_m(u),v))\otimes 1}}$$

$$\xymatrix{{V(\lambda(u,v))\otimes V(\mu(u,v))}\ar[r]^{1\otimes G_1(\mu(u,v))}&{V(\lambda(u,v))\otimes V(g_1(\mu(u,v)))}\ar[d]^{F(\lambda(u,v)),g_1(\mu(u,v)))}\\{V(u)\otimes V(v)}\ar[u]_{F(u,v)}\ar[d]^{1\otimes G_1(v)}&{V(\lambda(\lambda(u,v),g_1(\mu(u,v))))\otimes V(\mu(\lambda(u,v),g_1(\mu(u,v))))}\\{V(u)\otimes V(g_1(v))}\ar[r]^(.4){F(u,g_1(v))}&{V(\lambda(u,g_1(v)))\otimes V(\mu(u,g_1(v)))}\ar[u]_{1\otimes G_1(\mu(u,g_1(v)))}}$$

It is clear that we have also twisted transposition $\mu$ in this situation. By definition $\mu$ is the birational automorphism of $U\times U$ corresponding to the element $(1,m+1)(2,m+2)\dots(m-1,2m-1)\in S_{2m}$. We have 
$$\mu=\botshave\prod_{1\leqslant\alpha\leqslant m}g^l_{m-\alpha+1}\dots g^l_{m-1}fg^r_{m-1}\dots g^r_{\alpha}$$
where $g^l_{\alpha}=g_{\alpha}\times \text{id}$, $g^r_{\alpha}=\text{id}\times g_{\alpha}$ and the product is ordered from $\alpha=1$ to $\alpha=m$. One can also construct the corresponding twisted $R$-matrix as a product of $m^2$ factors:
$$R=\botshave\prod_{1\leqslant\alpha\leqslant m}(G_{m-\alpha+1}\otimes 1)\dots(G_{m-1}\otimes 1)F(1\otimes G_{m-1})\dots(1\otimes G_{\alpha}) \eqno(7)$$
Here we omit arguments of $R, G_{\alpha}$ and $F$.

Let $W$ be a vector space, $L$ be a meromorphic function from $U$ to $End(V\otimes W)$. We will consider $L(u)$ as a linear operator $L(u): V(u)\otimes W\to W\otimes V(u)$.

{\bf Definition} {\it We call $L$ an $L$-operator for twisted $GF$-matrices $G_i, F$ if the following diagrams are commutative:

$$\xymatrix{{V(u)\otimes W}\ar[rr]^{L(u)}\ar[d]^{G_i(u)\otimes 1}&&{W\otimes V(u)}\ar[d]^{1\otimes G_i(u)}\\{V(g_i(u))\otimes W}\ar[rr]^{L(g_i(u))}&&{W\otimes V(g_i(u))}}$$

$$\xymatrix{{V(u)\otimes W \otimes V(v)}\ar[rr]^{L(u)\otimes 1}&&{W\otimes V(u)\otimes V(v)}\ar[d]^{1\otimes F(u,v)}\\{V(u)\otimes V(v)\otimes W}\ar[u]_{1\otimes L(v)}\ar[d]^{F(u,v)\otimes 1}&&{W\otimes V(\lambda(u,v))\otimes V(\mu(u,v))}\\{V(\lambda(u,v))\otimes V(\mu(u,v))\otimes W}\ar[rr]^{1\otimes L(\mu(u,v))}&&{V(\lambda(u,v))\otimes W \otimes V(\mu(u,v))}\ar[u]_{L(\lambda(u,v))\otimes 1}}$$

That is}
$$G_i(u)L(u)=L(g_i(u))G_i(u), \text{ for } 1\leqslant i<m$$
$$F(u,v)L^1(u)L^2(v)=L^1(\lambda(u,v))L^2(\mu(u,v))F(u,v) \eqno(8)$$

We can consider $L(u)$ as a $n\times n$ matrix $L(u)=(L^j_i(u))$, where $L_i^j(u)$ are the elements in $End(W)$ such that 
$$L(u): x_i(u)\otimes W\to (L^j_i(u)W)\otimes x_j(u)$$
Here $\{x_i\}$ is a basis in $V$ and $\{x_i(u)\}$ is the corresponding basis in $V(u)$.

{\bf Definition} {\it The algebra of monodromy matrices $M_{GF}$ for the twisted $GF$-matrices $G_i, F$ is an associative algebra with the generators $\{L_i^j(u); i,j=1\dots n, u\in U\}$ and the defining relations (8).}

It is clear that the representations of the algebra $M_{GF}$ in the space $W$ coincide with the $L$-operators $L(u)\in End(V\otimes W)$.

Since the relations (4) follow from the relations (8), one has a homomorphism $M_{GF}\to M_R$, which is identity on the generators, where $R$ is given by (7). 

It is clear that the algebra $M_{GF}$ has the usual coproduct $\Delta: L^j_i(u)\to L_i^{\alpha}(u)\otimes L^j_{\alpha}(u)$.

{\bf Remarks 1.} One can also construct the Zamolodchikov algebra for $GF$-matrices.

{\bf 2.} For $Q(u)=\text{tr}L(u)=L_i^i(u)$ we have the following relations:
$$Q(u)=Q(g_i(u)) \text{ for } 1\leqslant i<m$$
$$Q(u)Q(v)=Q(\lambda(u,v))Q(\mu(u,v))$$

{\bf 3.} It is possible to construct twisted $GF$-matrices for the elliptic local action of the symmetric group defined in $\S4$ as intertwiners of tensor products of cyclic representations of the algebra of monodromy matrices for elliptic Belavin $R$-matrix ([5]) at the point of finite order. It will be the subject of another paper.

{\bf Acknowledgments}

I am grateful to V.Bazhanov and A.Belavin for useful discussions.

I am grateful to Max-Planck-Institut fur Mathematik, Bonn where this paper was written for invitation and very stimulating working atmosphere.

The work is supported partially by RFBR 99-01-01169, RFBR 00-15-96579, CRDF RP1-2254 and INTAS-00-00055.

\newpage
\centerline{\bf References}
\medskip

1. R.J. Baxter, J.H.H. Perk and H. Au-Yang, New solutions of the star-triangle relations for the chiral Potts model, Phys. Lett. A, 128(1988)138-142.

2. V.V. Bazhanov and Yu.G. Stroganov, Chiral Potts model as a descendant of the six vertex model, J. Stat. Phys. 51(1990)799-817.

3. V.O. Tarasov, Cyclic monodromy matrices for the $R$-matrix of the six-vertex model and the chiral Potts model with fixed spin boundary conditions, International Jornal of Modern Physics A, Vol.7, Suppl. 1B(1992)963-975.

4. A.A. Belavin, A.V. Odesskii and R.A. Usmanov, New relations in the algebra of the Baxter $Q$-operators. hep-th/0110126.

5. A.A. Belavin, Discrete groups and integrability of quantum systems, Func. Anal. Appl. 14(1980)18-26.

6. B.L. Feigin and A.V. Odesskii, Sklyanin Elliptic Algebras. The case of point of finite order, Func. Anal. Appl. 29(1995).

7. I. Gelfand, S. Gelfand, V. Retakh, S. Serconek, R.L. Wilson, Hilbert series of quadratic algebras associated with pseudo-roots of noncommutative polynomials, QA/0109007.

8. Mumford, Tata lectures on theta. I. With the assistance of C. Musili, M. Nory and M. Stillman. Progress in Mathematics, 28. BirkhDuser Boston, Inc., Boston, MA, 1983.

\enddocument